\documentclass{article}
\usepackage{amsmath,amssymb,color,authblk}
\usepackage[T1]{fontenc} 
\usepackage{pdfsync} 
\usepackage{latexsym}
\usepackage{amsmath}
\usepackage{amsfonts}
\usepackage{amsthm}

\newcommand{\hide}[1]{}

\definecolor{grey}{rgb}{0.75,0.75,0.75}
\definecolor{orange}{rgb}{1.0,0.5,0.5}
\definecolor{brown}{rgb}{0.5,0.25,0.0}
\definecolor{pink}{rgb}{1.0,0.5,0.5}

\definecolor{grey}{rgb}{0.75,0.75,0.75}
\definecolor{orange}{rgb}{1.0,0.5,0.5}
\definecolor{brown}{rgb}{0.5,0.25,0.0}
\definecolor{pink}{rgb}{1.0,0.5,0.5}
\newtheorem{theorem}{Theorem}
\newtheorem{lemma}[theorem]{Lemma}
\newtheorem{proposition}[theorem]{Proposition} 
\newtheorem{corollary}[theorem]{Corollary}

\newtheorem{notation}[theorem]{Notation} 
\newtheorem{definition}[theorem]{Definition}

\newtheorem{remark}[theorem]{Remark}

 \def\today{\number\day\space\ifcase\month\or
 janvier\or f\'evrier\or mars\or avril\or mai\or juin\or juillet\or
 ao\^ut\or septembre\or octobre\or novembre\or d\'ecembre\fi
 \space\number\year}

\def\ens#1{\left\{#1\right\}}
\def\Sylv{\mathrm{Sylv}}
\def\MSylv{\mathrm{MSylv}}
\def\BP{\mathbf{P}}
\def\BQ{\mathbf{Q}}

\def\BA{\mathbf{A} }
\def\BK{\mathbf{K} }
\def\BL{\mathbf{L} }
\def\BR{\mathbf{R} }
\def\BB{\mathbf{B} }

\def\BX{\mathbf{X} }
\def\BY{\mathbf{Y} }

\def\BT{\mathbf{T} }
\def\BU{\mathbf{U} }

\def\BW{\mathbf{W} }

\def\rem{\rm{Rem} }

\def\comb#1#2{\left(\begin{array}{c}#1\\#2 \end{array}\right)}

\def\Sres{{\rm{Sres}}}

\def\junk#1{{}}

\def\lc{\mathrm{lc}}

\def\C{{\rm{\bf C}}}

\newbox\auteurbox
\newbox\titrebox
\newbox\titrelbox
\newbox\editeurbox
\newbox\anneebox
\newbox\anneelbox
\newbox\journalbox
\newbox\volumebox
\newbox\pagesbox
\newbox\diversbox
\newbox\collectionbox
\def\fabriquebox#1#2{\par\egroup
\setbox#1=\vbox\bgroup
\leftskip=0pt
\hsize=\maxdimen \noindent#2}
\def\bibref#1{\bibitem{#1}
\mbox{}\ignorespaces
\setbox0=\vbox\bgroup}

{\catcode`\-=\active} 
{\catcode`\-=\active} 
{\catcode`\-=\active}
{\catcode`\-=\active}
\def\ajoutref#1{\setbox0=\vbox{\unvbox#1\global\setbox1=\lastbox}\unhbox1
\unskip\unskip\unpenalty}
\newif\ifpreviousitem
\global\previousitemfalse
\def\separateur{\ifpreviousitem {,\ }\fi} 
\def\voidallboxes {\setbox0=\box\auteurbox
\setbox0=\box\titrebox
\setbox0=\box\titrelbox
\setbox0=\box\editeurbox
\setbox0=\box\anneebox
\setbox0=\box\anneelbox
\setbox0=\box\journalbox
\setbox0=\box\volumebox
\setbox0=\box\pagesbox
\setbox0=\box\diversbox
\setbox0=\box\collectionbox
\setbox0=\null}
\def\fabriquelivre {\ifdim\ht\auteurbox>0pt
\ajoutref\auteurbox\global\previousitemtrue\fi 
\ifdim\ht\titrelbox>0pt
\separateur\ajoutref\titrelbox\global\previousitemtrue\fi 
\ifdim\ht\collectionbox>0pt
\separateur\ajoutref\collectionbox\global\previousitemtrue\fi 
\ifdim\ht\editeurbox>0pt
\separateur\ajoutref\editeurbox\global\previousitemtrue\fi 
\ifdim\ht\anneelbox>0pt	\separateur
\ajoutref\anneelbox
\fi\global\previousitemfalse}
\def\fabriquearticle {\ifdim\ht\auteurbox>0pt	\ajoutref\auteurbox
\global\previousitemtrue\fi
\ifdim\ht\titrebox>0pt
\separateur\ajoutref\titrebox\global\previousitemtrue\fi 
\ifdim\ht\titrelbox>0pt \separateur{\rm in}\
\ajoutref\titrelbox\global\previousitemtrue\fi 
\ifdim\ht\journalbox>0pt \separateur
\ajoutref\journalbox\global\previousitemtrue\fi 
\ifdim\ht\volumebox>0pt	\ \ajoutref\volumebox\fi
\ifdim\ht\anneebox>0pt	\ {\rm(}\ajoutref\anneebox \rm)\fi
\ifdim\ht\pagesbox>0pt
\separateur\ajoutref\pagesbox\fi\global\previousitemfalse}
\def\fabriquedivers {\ifdim\ht\auteurbox>0pt
\ajoutref\auteurbox\global\previousitemtrue\fi
\ifdim\ht\diversbox>0pt	\separateur\ajoutref\diversbox\fi}
\def\endbibref {\egroup
\ifdim\ht\journalbox>0pt \fabriquearticle 
\else\ifdim\ht\editeurbox>0pt \fabriquelivre 
\else\ifdim\ht\diversbox>0pt \fabriquedivers 
\fi\fi\fi .\voidallboxes}

\title{Sylvester double sums, subresultants and symmetric multivariate Hermite interpolation}
\date{}
\author[,1]{Marie-Fran\c{c}oise Roy}
\author[,2]{Aviva Szpirglas}
\affil[1]{Univ Rennes 1, CNRS, IRMAR-UMR 6625, F-35000 Rennes,  France}
\affil[2]{Univ Poitiers, CNRS, LMA-UMR 7348, F-86000 Poitiers, France}
\begin{document}
\maketitle

\begin{abstract}
Sylvester doubles sums, introduced first by Sylvester (see \cite{S1,S2}), are symmetric expressions of the roots of two polynomials $P$ and $Q$. Sylvester's definition of double sums makes no sense if $P$ an $Q$ have multiple roots, since the definition involves denominators that vanish when there are multiple roots. The aims of this paper are to give a new definition for Sylvester double sums making sense if $P$ and $ Q$ have multiple roots, which coincides with the definition by Sylvester in the case of simple roots, to prove the fundamental property of Sylvester double sums, i.e. that Sylvester double sums indexed by $(k,\ell)$ are equal up to a constant if they share the same value for $k+\ell$, and to prove the relationship between double sums and subresultants, i.e. that they are equal up to a constant.  In the simple root case, proofs of these properties are already known (see \cite{ALPP,AHKS,RS}). The more general proofs given here are using generalized Vandermonde determinants and a new symmetric multivariate Hermite interpolation as well as an induction on the length of the remainder sequence of $P$ and $Q$.
\end{abstract}

{\bf Keywords}: subresultants, Sylvester double sums, multivariate Hermite interpolation, generalized Vandermonde determinants

\section*{Introduction}

The first aim of this paper is to provide a definition for Sylvester double sums making sense if $P$ and $ Q$ have  multiple roots, which is done using quotients of Vandermonde determinants involving variables, and substitutions.
When the structure of the multiplicities of the roots of $P$ and $Q$ is known, we obtain a direct expresion of the Sylvester double sums in terms of generalized Vandermonde determinants.

The second aim of the paper is to prove, in the general case, the fundamental property for Sylvester double sums, i.e. that Sylvester double sums indexed by $(k,\ell)$ are equal up to a constant if they share the same value for $k+\ell$.
In order to prove this fundamental property, it is convenient to define more general objects, the {\sl multi Sylvester double sums}. We introduce a new multivariate symmetric Hermite interpolation and  use it to study the properties of  multi Sylvester double sums.  The strategy then consists in  proving the fundamental  property for multi Sylvester double sums and obtaining the result for Sylvester double sums as a corollary by identifying coefficients.  

The third aim of the paper is to prove the relationship between double sums and subresultants, i.e. that they are equal up to a constant. 
Our strategy is based on an induction on the length of the remainder sequence of $P$ and $Q$.

Our more general proofs are new even in the special case when the roots of the polynomials are simple.

The idea of introducing a multivariate symmetric Hermite interpolation and using multi Sylvester double sums  was inspired by \cite{KSV}'s use of multivariate symmetric Lagrange interpolation and introduction of multi Sysvtester double sums in the context of simple roots.

The content of the paper is the following.

In Section \ref{sec:def} we give a general definition for Sylvester double sums, valid also when there are multiple roots,
and prove that it coincides   with Sylvester's definition
in  the special case where all roots are simple (Proposition \ref{theoreme0}).

In Section \ref{sec:DSandVandermonde} we consider generalized Vandermonde determinants and use them to give a new formula for Sylvester double sums when the structure of multiplicities is known (Proposition \ref{lemme2bis}).

In Section \ref{sec:DSandHermite}, we introduce an Hermite interpolation for  multivariate symmetric polynomials (Proposition \ref{interpolationdata}).

In Section 
 \ref{sec:fonda} we study multi Sylvester double sums.
 We give their definition in subsection  \ref{subsec:defmulti}
In subsection \ref{subsec:mult} we compute the   multi Sylvester double sums and Sylvester double sums for indices $(k,\ell)$ with $k+\ell\ge \deg(Q)$.
In subsection \ref{subsec:fonda} we prove that multi Sylvester double sums and Sylvester double sums indexed by $k,\ell$, depend only (up to a constant) on $j=k+\ell$ (Theorem \ref{theo4mult} and Theorem \ref{theo4}). 

In Section \ref{sec:remainders} we give a relationship between Sylvester double sums of $(P,Q)$ and Sylvester double sums of  $(Q,R)$ where $R$ is the opposite of the remainder of $P$ by $Q$ in the Euclidean division (Proposition \ref{prorecurrence}).

 Finally we prove in Section \ref{sec:subresdblsum} that Sylvester double sums coincide (up to a constant) with subresultants, by an induction on the length of the remainder sequence of $P$ and $Q$ (Theorem \ref{theoreme2}).

\section{Sylvester double sums}  \label{sec:def}

We give a general definition for Sylvester double sums, valid also when the polynomials have multiple roots,
and prove that it coincides   with Sylvester's definition
in  the special case where all roots are simple (Proposition \ref{theoreme0}).

\subsection{Basic notations and definitions}

Let $\mathbb{K}$ be a field of characteristic $0$.

Let $\BA$
be
a finite 
list of elements of $\mathbb{K}$. 

We denote $\BA' \subset_a \BA$ when $\BA'$ is 
a sublist  of $\BA$ with $a$ elements (i.e. the list $\BA'$ is ordered by the restriction of the order on the list $\BA$). 

Let  $\BB$
be
another finite 
list of elements of $\mathbb{K}$.

We denote
$$\Pi(\BA,\BB)=\prod_{\genfrac{}{}{0pt}{}{x\in \BA}{y\in \BB}}(x-y).$$
Note that $\Pi(\BA, \BB)$ is independant on the order of $\BA$ and $\BB$.

We abbreviate  $\Pi(\{x\},\BB)$ and $\Pi(\BA,\{y\})$ to $\Pi(x,\BB)$ and $\Pi(\BA,y)$ respectively.

Note that $\Pi(\BA,\BB)$ is the classical resultant of the monic polynomials $\Pi(X,\BA)$ and $\Pi(X,\BB)$.

\begin{definition}\label{vderm}
The  {\rm Vandermonde vector }of length $i$ of $x\in \mathbb{K}$, denoted by $v_i(x)$, is
\begin{equation}
v_i(x)=\left[
\begin{array}{c}
1\\
x\\
\vdots\\
\vdots\\
x^{i-1}\\
\end{array}
\right].
\end{equation}
\end{definition}

\noindent
Let $\BA=(x_1,\ldots,x_i)$ be a finite ordered list of elements of $\mathbb{K}$.
The Vandermonde matrix $\mathcal{V}(\BA)$ is the $i \times i$ matrix having as  column vectors
$v_i(x_1),\ldots, v_i(x_i)$.
The Vandermonde determinant
$ V(\BA)$ is the determinant of the Vandermonde matrix  $\mathcal{V}(\BA)$.
It is well known that
\[ V(\BA)= \prod_{i \ge k > j \ge 1} (x_k - x_j) . \]

\medskip
By ${\bf B} \|{\bf A}$ the 
 we denote the
list obtained by concatening ${\bf B}$ and ${\bf A}$.

The following result is obvious.
\begin{lemma}\label{tresutile}
\begin{equation}
 V({\bf B} \|{\bf A})=V({\bf A}) \Pi ({\bf A},{\bf B}) V({\bf B}).
 \end{equation}
 and, as a special case, given a variable $U$,
   \[V({\bf B} \|U)=\Pi (U,{\bf B})V ({\bf B}).\]
\end{lemma}

\subsection{Definition of Sylvester double sums}

Let $\BP=(x_1,\ldots,x_p)$ and $\BQ=(y_1,\ldots,y_q)$ be  two finite ordered sets of element of $\mathbb{K}$
and $P=\Pi(X,\BP)$, $Q=\Pi(X,\BQ)$

The Sylvester double sum of $(P,Q)$ of index $k \in\mathbb{N}, \ell \in \mathbb{N} $
is usually defined 
as the following polynomial in $\mathbb{K}[U]$:

\begin{equation}\label{simple}
\sum_{\genfrac{}{}{0pt}{}{\BK \subset_k \BP}
      {\BL \subset_\ell \BQ}} \Pi (U,\BK) \Pi (U,\BL)
      \frac{\Pi (\BK, \BL) \Pi (\BP \setminus \BK,
      \BQ \setminus \BL)}{\Pi (\BK, \BP \setminus \BK)
      \Pi (\BL, \BQ \setminus \BL)} 
      \end{equation}
 (see \cite{S1,S2}).

This  definition of Sylvester double sums makes no sense if $P$ and $ Q$ have multiple roots, since
some of the quantities $\Pi (\BK, \BP \setminus \BK)$ (resp. $ \Pi (\BL, \BQ \setminus \BL)$) at the denominator are equal to $0$.

In this section, we give a general definition 
of Sylvester double sums, valid even if $P$ and $ Q$ have multiple roots and prove that it coincides with the classical one when all these roots are simple.

Let $\BX=(X_1,\ldots,X_p)$ and $\BY=(Y_1,\ldots,Y_q)$ be two ordered sets of indeterminates.

Given $\BX'\subset_k\BX$ (resp. $\BY'\subset_\ell\BY$ ), we denote $s_{\BX'}$  (resp.  $s_{\BY'}$ ) the signature of the permutation $\sigma_{\BX'}$ (resp. $\sigma_{\BY'}$)  putting
the elements of $\BX$ (resp. $\BY$) in the order $(\BX\setminus\BX')\|\BX'$ (resp.$(\BY\setminus\BY')\|\BY'$).

For any $k \in\mathbb{N}, \ell \in\mathbb{N}$,
we define the polynomial $F^{k,\ell}(\BX,\BY)(U)$ in $K[\BX,\BY,U]$
 \begin{equation}\label{defA} F^{k,\ell}(\BX,\BY)(U) =
     \sum_{\genfrac{}{}{0pt}{}{\BX' \subset_k \BX}
      {\BY' \subset_\ell \BY}} s_{\BX'} s_{\BY'}  V ((\BY \setminus \BY') \|
     (\BX\setminus \BX'))V (\BY' \| \BX' \| U)
     \end{equation}

Note that if $k>p$ or $\ell>q$ then $F^{k,\ell}(\BX,\BY)(U)=0$.
\begin{proposition}\label{lemme1}
The polynomial $F^{k,\ell}(\BX,\BY)(U)$ is antisymmetric in the variables $\BX$ and in the variables $\BY$.
\end{proposition} 

\begin{proof}
For any permutation $\sigma$ of the ordered set $\BX$, we call also $\sigma$
the action of $\sigma$ on a polynomial $F$ in $K[\BX,\BY,U]$, i.e $\sigma(F)(\BX,\BY)=F(\sigma(\BX),\BY)$. Denoting $s$ the signature
of $\sigma$ we want to prove 
\begin{equation}\label{anti}
\sigma(F^{k,\ell})(\BX,\BY)(U)=s F^{k,\ell}(\BX,\BY)(U).
\end{equation}
It is enough to prove (\ref{anti}) for a transposition exchanging two sucessive elements, of  signature $-1$.

So, let  $\tau$ be the transposition exchanging $X_i$ and $X_{i+1}$.
We want to prove
\begin{equation}\label{anti1}
\tau(F^{k,\ell})(\BX,\BY)(U)=-F^{k,\ell}(\BX,\BY)(U).
\end{equation}

 We denote by $\tau(\BX)$ the ordered set obtained from
$\BX$ by exchanging $X_i$ and $X_{i+1}$. Given $\BX'\subset_k\BX$,
we denote  by $\tau(\BX')$ 
the ordered set  $\tau(\BX')\subset_k \tau(\BX)$ (i.e. $\tau(\BX')$ is ordered by the restriction of the order on $\tau(\BX)$)
and by $\bar \BX'$ the ordered set $\tau
(\BX')\subset_k \BX$ (i.e. $\bar \BX'$ and $\tau
(\BX')$ have the same elements but  $\bar \BX'$ is ordered by the restriction of the order on $\BX$).

Denote
$$
F^{\BX',\BY'}=s_{\BX'} s_{\BY'}V ((\BY \setminus \BY') \|
     (\BX\setminus \BX'))V (\BY' \| \BX' \|U) .
$$

We have 3 cases to consider.

\begin{itemize}
\item If $X_i\in \BX'$  and $X_{i+1}\in \BX'$ then $\tau(\BX \setminus \BX')= \BX \setminus \BX'$ and
$$\begin{array}{rcl}
\tau(F^{\BX',\BY'})&=&s_{\BX'} s_{\BY'}
V ((\BY \setminus  \BY') \|
    \tau(\BX \setminus \BX'))V ( \BY'\| \tau(\BX') \|U) \\
                                    &=&s_{\BX'} s_{\BY'}V ((\BY \setminus \BY') \|
    ( \BX \setminus \BX'))V ( \BY' \| \tau(\BX') \|U) \\
                                     &=&-F^{\BX',\BY'}.
\end{array}$$
\item If $X_i\notin \BX'$  and $X_{i+1}\notin \BX'$ then $\tau(\BX')=\BX'$ and
$$\begin{array}{rcl}
\tau(F^{\BX',\BY'})&=&s_{\BX'} s_{\BY'}
V ((\BY \setminus  \BY') \|
    \tau(\BX \setminus \BX'))V ( \BY'\| \tau(\BX') \|U) \\
                                    &=&s_{\BX'} s_{\BY'}V ((\BY \setminus \BY') \|
    \tau( \BX \setminus \BX'))V ( \BY' \| \BX' \|U) \\
                                     &=&-F^{\BX',\BY'}.
\end{array}$$
\item If $X_i \in \BX'$ and $X_{i+1} \notin \BX'$ , or  $X_i \notin \BX'$ and $X_{i+1} \in \BX'$, then  
  $\sigma_{\bar \BX'}=\tau \circ \sigma_{\BX'}$, $\tau(\BX')=\bar \BX'$ and  $\tau(\BX \setminus \BX')=\BX \setminus \bar \BX'$
     so that
$$\begin{array}{rcl}
\tau(F^{\BX',\BY'})&=&s_{\BX'} s_{\BY'} V ( (\BY \setminus  \BY') \|
    \tau(\BX\setminus \BX'))V (\BY' \| \tau(\BX') \|U)\\
                                   &=&- s_{\bar \BX'} s_{\BY'} V ((\BY\setminus \BY') \|
     (\BX \setminus \bar \BX'))V (\BY' \| \bar \BX' \|U)\\
                                   &=&-F^{\bar \BX',\BY'}
\end{array}$$
and 
 $$\begin{array}{rcl}
\tau(F^{\bar \BX',\BY'})&=&s_{\bar \BX} s_{\BY'}V (( \BY \setminus \BY') \|
 \tau (  \BX \setminus \bar\BX'))V ( \BY'\|\tau(\bar\BX') \|U) \\
                                     &=&-s_{\BX} s_{\BY'}V ((\BY \setminus \BY') \|
     (\BX \setminus \BX'))V (\BY' \| \BX' \|U) \\
                                     &=&-F^{\BX',\BY'},
\end{array}$$
From which we deduce
$$
\tau
\left(F^{\BX',\BY'}+F^{\bar \BX',\BY'}\right)=-\left(F^{\bar \BX',\BY'}+F^{\BX',\BY'}\right).
$$

So, we get (\ref{anti1}). 

\end{itemize}

The exchange between two elements of $\BY$ can be treated similarly.
\end{proof}

\begin{lemma}\label{div}
  If $A(\BX,\BY)$ in $K[\BX,\BY]$ is antisymmetric with respect to the variables $\BX$, then $A(\BX,\BY)=S(\BX,\BY) V(\BX)$ where $S\in K[\BX,\BY]$ is symmetric with respect to the variables $\BX$.
\end{lemma}
\begin{proof}
If $A(\BX,\BY)$ is antisymmetric with respect to $\BX$ then, for any $j<k$, 
denote  $\tau_{j,k}(\BX)$ the ordered set of variables obtained by transposing $X_j$ and $X_k$.

\medskip
$$
\frac{A(\BX,\BY)-A(\tau_{j,k}(\BX,\BY))}{X_j-X_k}=2\frac{A(\BX,\BY)}{X_j-X_k}
$$
 is a polynomial.
So $A(\BX,\BY)=S(\BX,\BY) V(\BX) $  and $S(\BX,\BY)$ is a symmetric polynomial with respect to $\BX$.
\end{proof}

Applying Lemma \ref{div} and Proposition \ref{lemme1} we denote
 $S^{k,\ell}(\BX,\BY)(U)$  the symmetric polynomial with respect to the indeterminates $\BX$ and with respect to the indeterminates $\BY$ satisfying 
\begin{equation}\label{defS}
S^{k,\ell}(\BX,\BY)(U)=\frac{F^{k,\ell}(\BX,\BY)(U)}{V(\BX)V(\BY)}.
\end{equation}

Given two monic univariate polynomials $P$ and $Q$ of degree $p$ and $q$ we denote 
$\BP=(x_1,\ldots,x_p)$ and $\BQ=(y_1,\ldots,y_q)$ ordered lists of the roots of $P$ and $Q$ in an algebraic closure $\C$ of $\mathbb{K}$, counted with multiplicities.

\begin{definition} The  {\rm generalized Sylvester double sum of $(P,Q)$} for the exponents $k,\ell\in\mathbb{N}\times\mathbb{N}$ is defined by 
$$\Sylv^{k,\ell}(P,Q)(U)=S^{k,\ell}(\BP,\BQ)(U).$$
\end{definition}
Note that this definition does not depend on the order given for the roots of  $P$ and  $Q$.

\medskip This definition of generalized Sylvester double sums for monic polynomials coincides  with the usual definition of Sylvester double sums  when the  polynomials $P$ 
and $Q$ have no multiple roots, as we see now in Proposition \ref{theoreme0}.

\noindent 
\begin{proposition}\label{theoreme0} If $P,Q$ have only simple roots, 
$$\Sylv^{k,\ell}(P,Q )(U)=\sum_{\genfrac{}{}{0pt}{}{\BK \subset_k \BP}
      {\BL \subset_\ell \BQ}} \Pi (U,\BK) \Pi (U,\BL)
      \frac{\Pi (\BK, \BL) \Pi (\BP \setminus \BK,
      \BQ \setminus \BL)}{\Pi (\BK, \BP \setminus \BK)
      \Pi (\BL, \BQ \setminus \BL)}$$
\end{proposition}

\break
\begin{proof}[Proof of Proposition \ref{theoreme0}]\

\noindent $\displaystyle{ \sum_{\genfrac{}{}{0pt}{}{\BK \subset_k \BP}
      {\BL \subset_\ell \BQ}} \Pi (U,\BK) \Pi (U,\BL)
      \frac{\Pi (\BK, \BL) \Pi (\BP \setminus \BK,
      \BQ \setminus \BL)}{\Pi (\BK,\BP \setminus \BK)
      \Pi (\BL, \BQ \setminus \BL)}}=$
$$\begin{array}{cl}
=&  \displaystyle{\sum_{\genfrac{}{}{0pt}{}{\BK \subset_k \BP}
      {\BL \subset_\ell \BQ}} \Pi (U,\BK) \Pi (U,\BL)
      \frac{V(\BL) \Pi (\BK, \BL) V(\BK) V(\BQ \setminus \BL) \Pi (\BP \setminus \BK,\BQ \setminus \BL)V(\BP \setminus \BK)}{V(\BK)\Pi (\BK,\BP \setminus \BK)  V(\BP \setminus \BK)
      V(\BL) \Pi (\BL,\BQ \setminus \BL) V(\BQ \setminus \BL) }}\\
  =&    \displaystyle{ \sum_{\genfrac{}{}{0pt}{}{\BK \subset_k \BP}
      {\BL \subset_\ell \BQ}}s_{\BK} s_{\BL}  \frac{V (\BL  \| \BK  \|
   U) V ((\BQ \setminus \BL)\| (\BP \setminus \BK))}{V (\BP)
   V (\BQ)}}\\
 =&  \displaystyle{\frac{F^{k,\ell}(\BP,\BQ)(U)}{V(\BP)V(\BQ)}}=S^{k,\ell}(\BP,\BQ)(U)\\\\
 =&\Sylv^{k,\ell}(P,Q)(U)
\end{array}$$
applying Lemma \ref{tresutile}.
 \end{proof}

\section{Generalized Vandermonde determinants
and Sylvester double sums}\label{sec:DSandVandermonde}

We consider generalized Vandermonde determinants (also called sometimes confluent Vandermonde determinants, see \cite{LT,HJ})  and connect them with the Sylvester double sums (Proposition \ref{lemme2bis}).

\begin{notation}\label{notmultiset}
 Let  $P$ be a 
 polynomial of degree $p$ with coefficients in a field $\mathbb{K}$. Let $(x_1,\ldots,x_m)$  be an ordered set of the distinct roots of $P$ in an algebraic closure $\C$ of $\mathbb{K}$, with $x_i$ of multiplicity $\mu_i$, and let $\BP$ be the  multiset of roots of $P$, represented by the ordered set
$$\BP=(x_{1,0},\ldots,x_{1,\mu_1-1},\ldots,x_{m,0},\ldots,x_{m,\mu_m-1}), $$
with $x_{i,j}=(x_i,j)$ for $0\leq j\leq \mu_i-1$, $\sum_{i=1}^m \mu_i=p.$

\medskip Let  $Q$ be a 
polynomial of degree $q$ with coefficients in $\mathbb{K}$. Let $(y_1,\ldots,y_n)$  be an ordered set of the distinct roots of $Q$ in $\C$
 with $y_i$ of multiplicity $\nu_i$, , for $i=1,\ldots,n$. Let  $\BQ$ be  the ordered multiset of its root, represented by the ordered set
$$\BQ=(y_{1,0},\ldots,y_{1,\nu_1-1},\ldots,y_{n,0},\ldots,y_{n,\nu_n-1}), $$
with $y_{i,j}=(y_i,j)$ for $0\leq j\leq \nu_i-1$,  $\sum_{i=1}^n \nu_i=q.$

 We introduce an ordered set of variables $X_\BP=(X_{1,0}, \ldots,X_{1,\mu_1-1},\ldots,X_{m,0},\ldots,X_{m,\mu_m-1})$ and an ordereed set of variables $Y_\BQ=(Y_{1,0}, \ldots,Y_{1,\nu_1-1},\ldots,Y_{n,0},\ldots,Y_{n,\nu_n-1})$.
 
 For a polynomial $f(X_\BP,X_\BQ)$ we denote $f(\BP,\BQ)$ the result of the substitution of $X_{i,j}$ by $x_i$ and $Y_{i,j}$ by $y_i$.

\end{notation}

\begin{notation}
Given $f$  a  polynomial depending on the variable $U$, we denote 
\begin{equation}f^{[i]}=\frac1{i!}\frac{\partial^{i} f}{\partial U^i }.\label{derdermult}
\end{equation}
\end{notation}

\begin{definition}[Generalized Vandermonde determinant]
Let  $\BK\subset_k \BP$, $\BL \subset_\ell \BQ$ and $\BU=(U_1,\ldots,U_u)$ an ordered set of $u$ indeterminates.

\medskip The {\rm generalized Vandermonde matrix}
${\cal V}[\BL\|\BK\|\BU)$ is the $(\ell+k+u) \times (\ell+k+u)$ matrix having as column vectors   the $\ell$ columns
$v^{[j]}_{k+\ell+u}(y_i)$ for $y_{i,j}\in \BL$
followed by the $k$ columns
$v^{[j]}_{k+\ell+u}(x_i)$ for $x_{i,j}\in \BK$ 
followed by the $u$ columns $v_{k+\ell+u}(U_i)$
(using notation (\ref{vderm}) and notation (\ref{derdermult})). 

\medskip The {\rm generalized Vandermonde determinant} $V[\BL\|\BK\|\BU)$ is the determinant of ${\cal V}[\BL\|\BK|\BU)$.

\begin{itemize}
\item In the particular case $u=0$ we denote $V[\BL\|\BK]$ the corresponding determinant.

\item In the particular case $k=p,\ell=u=0$ we denote $V[\BP]$ the corresponding determinant .

\item Similarly, in  the particular case $k=0,\ell=q,u=0$ we denote $V[\BQ]$ the corresponding determinant.
\end{itemize}

\begin{remark}{}
{\rm The peculiar notation $V[\BL\|\BK\|\BU)$ with one square bracket to the left and one parenthesis to the right is here to indicate that the column $v_{k+\ell+u}^{[j]}(x_i)$ indexed by $x_{i,j}\in \BK$ and  $v_{k+\ell+u}^{[j]}(y_i)$ indexed by $y_{i,j}\in\BL$ have been derivated, while there is no derivation with respect to the columns indexed by the variables in $\BU$.}
\end{remark}
\end{definition}

While the classical Vandermonde determinant $V(\BP)$ is null when $P$ has multiple roots, we have the following result for the generalized 
Vandermonde determinant.

\begin{lemma}
 The  generalized Vandermonde determinant $ V[\BP]$ is equal to
$$
 V[\BP]
=\prod_{1\leq i<j\leq m }(x_j-x_i)^{\mu_i\mu_j}.$$ 
\end{lemma}
\begin{proof}
The proof is done by induction on $p$.

If $p=1$, $V[\BP]=1$.

Suppose that $$\displaystyle{ V[\BP]
=\prod_{1\leq i<j\leq m }(x_j-x_i)^{\mu_i\mu_j}}.$$

The polynomial $F(U)=V[\BP\| U)$  is  of degree $p$, with leading coefficient $V[\BP]$   and satisfies the property
 $$
\mbox{ for all  }1\leq i\leq m, \mbox{ for all  }0\leq j < \mu_i, 
\begin{array}{lcl}
F^{[j]}(x_i)&=&0,
\end{array}
$$
So, $$F(U)=V[\BP]\prod_{i=1}^{m}(U-x_i)^{\mu_i}=\prod_{1\leq i<j\leq m }(x_j-x_i)^{\mu_i\mu_j}\prod_{i=1}^{m}(U-x_i)^{\mu_i}.$$

 Consider $T(U)=(U-x)P(U)$. 
 
\medskip
$-$ First case: $x$ is not a root of $P$. Let  $\BT$ the ordered set (obtained by adding $x$ at the end of $\BP$) of roots of the polynomial $T$, so $x=x_{m+1}$ is a root of $T$ with multiplicity 1.
Then $$V[\BT]=F(x)=\prod_{1\leq i<j\leq m }(x_j-x_i)^{\mu_i\mu_j}\prod_{i=1}^{m}(x-x_i)^{\mu_i}=\prod_{1\leq i<j\leq m+1 }(x_j-x_i)^{\mu_i\mu_j}.$$

\medskip
$-$ Second case: $x$ is  a root of $P$. So there exists
$1\leq j\leq m$ such that $x=x_j$ , and $x_j$ is a root of multiplicity $\mu_j+1$ of $T$.

Let  $\BT$ the ordered set of roots of the polynomial $T$ obtained by inserting $x_{j,\mu_j}=(x_j,\mu_j)$ after
 $x_{j,\mu_j-1}$ in $\BP$.
 Then
\begin{align*}
V[\BT]&=(-1)^{\mu_{j+1}+\cdots+\mu_m}F^{[\mu_j]}(x_j)\\
&=(-1)^{\mu_{j+1}+\cdots+\mu_m}V[\BP]\displaystyle{\prod_{\genfrac{}{}{0pt}{}{i=1}{i\not= j}}^{m}(x_j-x_i)^{\mu_i}}\\
&=\displaystyle{\prod_{1\leq i<j\leq m }(x_j-x_i)^{\mu_i(\mu_j+1)}}\qedhere
\end{align*}
\qedhere
\end{proof}

\begin{remark}{}
{\rm If $\BK\subset_k \BP$, it can happen that $V[\BK]=0$.
 Taking for example $\BP=(x_{1,0},x_{1,1},x_{2,0},x_{2,1})$ and $\BK=(x_{1,1},x_{2,1})$,  it is easy to check that $ V[\BK]=0$.
}
\end{remark}

From now on, and till Section \ref{sec:remainders},
$P$ and $Q$ are monic polynomials, 
 
 \medskip The following proposition makes the link between generalized Vandermonde determinants and Sylvester double sums.

We denote  $s_{\BK}$ (resp.  $s_{\BL}$)  the signature of the permutation $\sigma_{\BK}$ 
(resp.  $\sigma_{\BL}$) obtained by putting
the elements of $\mathbf{P}$ (resp.  $\BQ$) in the order $(\BP\setminus\BK)\|\BK$ (resp. $(\BQ\setminus\BL)\|\BL$).

\begin{proposition}\label{lemme2bis} ~
$$
\Sylv^{k,\ell}(P,Q)(U)= \sum_{\genfrac{}{}{0pt}{}{\BK \subset_k \BP}
      {\BL \subset_\ell \BQ}} s_{\BK} s_{\BL}\frac{
  V[(\BQ \setminus \BL) \|
     (\BP \setminus \BK)] V[\BL \| \BK \| U) } {V[\BP]  V[\BQ]}
$$
\end{proposition}
In the proof of Proposition \ref{lemme2bis}, we use the following notation \ref{derivations}
and Lemma \ref{exemple2}.

\begin{notation}
\label{derivations}
For any polynomial $f$ depending on the variables $X_\BP$, and $\BK \subset_k \BP$, denote $\partial^{[\BK]}f$ 
the polynomial defined by induction on $r$
as follows. 
$$\partial^{[\emptyset]}f=f$$
If $\BK=\BK'\|(x_{i,j})$, 
$$\partial^{[\BK]}f=\frac1{j!}\frac{\partial^{j}\partial^{[\BK']}f}{\partial X_{i,j}^j }.$$

\medskip Similarly, for any polynomial $f$ depending on the variables $Y_\BQ$, and $\BL \subset_\ell \BQ$, denote $\partial^{[\BL]}f$ 
the polynomial defined by induction on $s$ as follows. 
$$\partial^{[\emptyset]}f=f$$
If $\BL=\BL'\|(y_{i,j})$, 
$$\partial^{[\BL]}f=\frac1{j!}\frac{\partial^{j}\partial^{[\BL']}f}{\partial X_{i,j}^j }.$$

Note that 
$$V[\BL\|\BK\|\BU)=f(\BK,\BL,\BU)$$
with $f(X_\BK,Y_\BL,\BU)=\partial^{[\BL]}\partial^{[\BK]}V(Y_\BL\|X_\BK\|\BU)$.

\end{notation}

\begin{lemma}\label{exemple2}~
$$\partial^{[\BP]}\left(V(X_\BP)f(X_\BP)\right)(\BP)= V[\BP] f(\BP)$$
\end{lemma}
\begin{proof}
We first note that 
$$\partial^{[\BP]}\left(V(X_\BP)f(X_\BP)\right)=\partial^{[\BP]}(V(X_\BP))f(X_\BP)+\sum_r V_r(X_\BP) f_r(X_\BP)$$
where   $V_r(X_\BP)$ (resp. $f_r(X_{\BP}^{})$) is obtained from $V(X_\BP)$ (resp. from $f(X_\BP)$) by partial derivations, one variable $X_{i,j}$ at least being derived less than $j$ times (resp. at least one time). Denoting 
$X_{i,j}$ the first variable which is being derived less than $j$ times in $V_r(X_\BP)$, we define $j'$ as the order of derivation of $X_{i,j}$ in  $V_r(X_\BP)$. We notice that 
$V_r(X_\BP)$ is the determinant of a matrix with two equal columns, the one indexed by ${i,j'}$ and 
the one indexed by ${i,j}$. Hence $V_r(\BP)=0$. This proves the claim.
\end{proof}

\begin{proof} [Proof of Proposition \ref{lemme2bis}]

\medskip Since
$$F^{k,\ell}(X_\BP,Y_\BQ)(U)=V(X_\BP)V(Y_\BQ)S^{k,\ell}(X_\BP,Y_\BQ)(U),$$ using Lemma \ref{exemple2} we obtain
$$\partial^{[\BQ]}\partial^{[\BP]}F^{k,\ell}(\BP,\BQ)(U)= V[\BP] V[\BQ]S^{k,\ell}(\BP,\BQ)(U)= V[\BP] V[\BQ]\Sylv^{k,\ell}(P,Q)(U).$$
On the other hand, 
denoting
$$h_{\BK,\BL}(X_\BP,Y_\BQ)(U)= 
 V (Y_{\BQ \setminus \BL} \|
     X_{\BP \setminus\BK})V (Y_\BL\| X_\BK \| U),$$ we have
     $$\partial^{[\BQ]}\partial^{[\BP]} h_{\BK,\BL}(\BP,\BQ)(U)= V[(\BQ \setminus \BL) \| (\BP \setminus \BK)]V[\BL \| \BK \| U).$$
   
Since  
$$
F^{k,\ell}(X_\BP,Y_\BQ)(U)=\sum_{\genfrac{}{}{0pt}{}{\BK \subset_k \BP}{\BL \subset_\ell \BQ}}s_\BK s_\BL h_{\BK,\BL}(X_\BP,Y_\BQ)(U)
$$
we get
\begin{align*}
\partial^{[\BQ]}\partial^{[\BP]}F^{k,\ell}(\BP,\BQ)(U)&=\sum_{\genfrac{}{}{0pt}{}{\BK \subset_k \BP}
      {\BL \subset_\ell \BQ}} s_{\BK} s_{\BL}   V[(\BQ \setminus \BL )\| (\BP \setminus \BK)]V[\BL \| \BK \| U)\qedhere
\end{align*}
\end{proof}

The following lemma will be useful later.

\begin{lemma}\label{prodpratique} ~
\begin{enumerate}
\item  For $\BL\subset_\ell\BQ$,
  defining 
$$f(Y_{\BQ\setminus \BL})=(-1)^{p(q-\ell)}\partial^{{[\BQ\setminus\BL]}}\left(V(Y_{\BQ\setminus\BL})\prod_{Y\in Y_{\BQ\setminus\BL}}P(Y)\right),$$
we have
$$
 V[(\BQ \setminus \BL)\|\BP ]=V[\BP]f(\BQ\setminus\BL).
$$
\item For $\BK\subset_k\BP$,
 defining
$$g( X_{\BP\setminus\BK})=\partial^{{[\BP\setminus\BK]}}\left(V(X_{\BP\setminus\BK})\prod_{X \in X_{\BP\setminus\BK}}Q(X) \right),$$ 
we have
$$
 V[\BQ\|(\BP\setminus\BK)]=V[\BQ]g{(\BP\setminus\BK)}
$$
\end{enumerate}

\end{lemma}
\begin{proof}
Defining 
$$\begin{array}{rcl}
h(X_\BP,Y_{\BQ\setminus\BL})&=&\partial^{[\BQ\setminus\BL]}V(Y_{\BQ\setminus\BL}\|X_\BP)
\\
&=&\partial^{[{\BQ\setminus\BL}]}\left(V(X_\BP)\Pi(X_\BP,Y_{\BQ\setminus\BL})V(Y_{\BQ\setminus\BL})\right)\\
&=&V(X_\BP)\partial^{[{\BQ\setminus\BL}]}\left(\Pi(X_\BP,Y_{\BQ\setminus\BL}) V(Y_{\BQ\setminus\BL})\right)
\end{array}$$
and applying Lemma \ref{exemple2},
we get
$$\begin{array}{rcl}
\partial^{[\BP]}h(\BP,Y_{\BQ\setminus\BL})&=&V[\BP]\partial^{[{\BQ\setminus\BL}]}\left(V(Y_{\BQ\setminus\BL})\Pi(\BP,Y_{\BQ\setminus\BL})\right)\\
&=&V[\BP]\partial^{[{\BQ\setminus\BL}]}\left(V(Y_{\BQ\setminus\BL})\displaystyle{\prod_{Y\in Y_{\BQ\setminus\BL}}(-1)^pP(Y)}\right)\\
&=&V[\BP]f(Y_{\BQ\setminus \BL}) .
\end{array}$$
and finally
$$V[(\BQ\setminus\BL) \|\BP]=\partial^{[\BP]}h(\BP,\BQ\setminus\BL)=V[\BP]f(\BQ\setminus\BL).
$$
Which is Lemma \ref{prodpratique}.1.

The proof for Lemma \ref{prodpratique}.2 is similar.
\end{proof}

Lemma \ref{prodpratique} has the following corollary.

\begin{corollary}\label{Regal0}
If 
$Q$ divides $P$, then $$\Sylv^{0,j}(P,Q)(U)=0$$
\end{corollary}
\begin{proof}
In this case, $f(\BQ\setminus\BL)=0$ as any root of $Q$ is a root of $P$ with at least the same multiplicity. So, applying Lemma \ref{prodpratique}.1, $V[(\BQ\setminus\BL)\|\BP]=0$.
It follows
$$
\Sylv^{0,j}(P,Q)(U)=\sum_{\BL\subset_j\BQ}^{}s_\BL\frac{V[(\BQ\setminus\BL)\|\BP]V[\BL\|U)}{V[\BP]V[\BQ]}=0$$
\end{proof}
\section{Hermite Interpolation for multivariate symetric polynomials}
\label{sec:DSandHermite}
\label{symmetrichermite}

We now introduce an Hermite interpolation for  multivariate symmetric polynomials.

\medskip We consider an ordered set of $p-k$ variables $\BU$.

\begin{proposition}\label{multih}
The set $$\mathcal{B}_{\BP,k}(\BU)=\ens{\frac{ V[\BK\|\BU)}{V[\BP]V(\BU)}\mid \BK\subset_k \BP }$$ 
is a  basis of the vector-space of symmetric polynomials in $\BU$ of multidegree at most $k,\ldots,k$.
\end{proposition}

The proof of Proposition  \ref{multih} uses the following Lemma.

\begin{lemma}\label{toutourien}\

\begin{enumerate}
\item $ V[\BK\|(\BP\setminus \BK)]=(-1)^{k(p-k)}s_\BK  V[\BP]\not=0.$
\item If $\BK'\not= \BK$, $V[\BK'\|(\BP\setminus \BK) ]=0.$
\end{enumerate}

\end{lemma}
\begin{proof}\

\begin{enumerate}
\item It is clear that
$V[\BK\|(\BP\setminus \BK)]=(-1)^{k(p-k)}s_\BK  V[\BP]\not=0,$
since $s_\BK$ is the signature of the permutation putting $\BP$ in the order $(\BP\setminus\BK)\| \BK$.

\item The fact that
$V[\BK'\|(\BP\setminus \BK)]=0$
when $\BK'\not= \BK$ follows from  the fact that  the matrix $ {\cal V}[\BK'\|(\BP\setminus \BK)]$ has two equal columns.\qedhere
\end{enumerate}

\end{proof}

\begin{proof} [Proof of Proposition  \ref{multih}]
Since the number of subsets of cardinality $k$ of $\BP$ is $\displaystyle{\binom{p}{k}}$ and that $\displaystyle{\binom{p}{k}}$ is also the dimension of the vector space of symmetric polynomials in $\BU$ of multidegree at most $k,\ldots,k$,  it is enough to prove that 
$$\sum_{\BK' \subset_k \BP} c_{\BK'}\frac{ V[\BK'\|\BU)}{V[\BP]V(\BU)}=0$$
implies $c_{\BK}=0$ for all $\BK\subset_k \BP$.

Let us fix $\BK\subset_k \BP$.
Since
$$\sum_{\BK' \subset_k \BP} c_{\BK'}  V[\BK'\|\BU)=0,$$
it follows by substitution and derivation that
$$\sum_{\BK' \subset_k \BP} c_{\BK'} \partial^{[\BP\setminus \BK]} V[\BK'\|X_{\BP\setminus \BK})=0.$$
When replacing $X_{\BP\setminus \BK}$ by $\BP\setminus \BK$ we obtain
\begin{align*}
\sum_{\BK' \subset_k \BP} c_{\BK'} V[\BK'\|(\BP\setminus \BK)]&=0.
\end{align*}
Using Lemma \ref{toutourien}, we get
$c_{\BK}=0$.
\end{proof}

The following Proposition gives the connection between a symetric polynomial in $\BU$ of multidegree at most $k,\ldots,k$ and its coordinates in the basis $\mathcal{B}_{\BP,k}(\BU)$.

\begin{proposition}{\bf (Multivariate symmetric Hermite Interpolation)} \label{interpolationdata}
Let $g$  be a symetric polynomial in $\BU$ of multidegree at most $k,\ldots,k$.
Writing
 $$g(\BU)=\sum_{\BK \subset_k \BP} g_{\BK}\frac{  V[\BK\|\BU)}{V[\BP]V(\BU)}$$
then
$$g_\BK=(-1)^{k(p-k)}s_\BK  h(\BP\setminus \BK)$$ 
with
$$h(X_{\BP\setminus \BK})= \partial^{[\BP\setminus \BK]}(V(X_{\BP\setminus \BK})g(X_{\BP
\setminus \BK}) ).$$
\end{proposition}
\begin{proof}
We have
$$\sum_{\BK \subset_k \BP} g_\BK   V[\BK\|\BU)=V[\BP]g(\BU)V(\BU).$$
Derivating both sides by $\partial^{[\BP\setminus \BK']}$ and substituting $\BP\setminus \BK'$ for $\BU$ we get, using Lemma \ref{toutourien}
$$g_{\BK' } V[\BK'\|(\BP\setminus \BK')]=g_{\BK' }s_{\BK'}  (-1)^{k(p-k)} V[\BP]= V[\BP]h(\BP\setminus \BK') ,$$
and finally
\[g_{\BK'}=(-1)^{k(p-k)}s_{\BK'} h(\BP\setminus \BK').\qedhere\]
\end{proof}

\begin{remark}{}
{\rm Proposition \ref{multih} generalizes a result in 
\cite{CL} 
given for Lagrange interpolation of symmetric multivariate polynomials.}
\end{remark}

As a corollary of Proposition \ref{interpolationdata}, we recover the classical Hermite Interpolation

\begin{proposition}\label{hermite}{\bf (Hermite Interpolation)}
Given an ordered list
 $${\bf q}=(q_{1,0},\ldots,q_{1,\mu_1-1},\ldots,q_{m,0},\ldots,q_{m,\mu_m-1})$$ of $p$ numbers,
there is one and  only one polynomial  of degree at most $p-1$ satisfying the property
$$
\mbox{ for all  }1\leq i\leq m, \mbox{ for all  }0\leq j < \mu_i, 
\begin{array}{lcl}
Q^{[j]}(x_i)&=&q_{i,j}.
\end{array}
$$
\end{proposition}

\begin{proof}
If $k=p-1$ in Proposition \ref{multih}, then 
$$
   \mathcal{B}_{\BP,p-1}(U)= \ens{\frac{V[\BP\setminus \ens{x_{i,j}}\|U)}{V[\BP]}\mid x_{i,j}\in\BP}
    $$
     is a basis of the vector space of univariate polynomials in $U$  of degree at most $\le p-1$.
 
Note that $(-1)^{p-1}s_{\BP\setminus \ens{x_{i,j}}}=(-1)^{\mu_{i}+\cdots+\mu_m-j-1}$. So, the family $$\ens{(-1)^{\mu_{i}+\cdots+\mu_m-j-1}q_{i,j}\mid i=1,\ldots,m, j=0,\ldots \mu_i-1}$$ is the coordinates in the basis $\mathcal{B}_{\BP,p-1}(U)$ of a polynomial $Q(U)$ (necessarily unique) of degree at most $p-1$ such that $Q^{[j]}(x_i)=q_{i,j}$, applying Proposition \ref{interpolationdata}. 
\end{proof}

\section{Multi Sylvester double sums} \label{sec:fonda}

We introduce in subsection 
 \ref{subsec:defmulti} multi Sylvester double sums and study their properties, using  the  Hermite interpolation for symmetric multivariate polynomials.
In subsection \ref{subsec:mult} we compute the   multi Sylvester double sums and Sylvester double sums for indices $(k,\ell)$ with $k+\ell\ge q$.
In subsection \ref{subsec:fonda} we prove the fundamental property of Sylvester double sums, i.e. that Sylvester double sums indexed by $k,\ell$, depend only (up to a constant) on $j=k+\ell<p$. This was already known in the simple roots case but even in this case our proof is new. 

\subsection{Definition of multi Sylvester double sums} 
\label{subsec:defmulti}

The idea of replacing the variable $U$ by a block of indeterminates to define  multi Sylvester double sums is directly inspired from  \cite{KSV}.

\begin{definition}
The {\rm multi Sylvester double sum}, for $(k,\ell)$ a pair of natural numbers with $ k+\ell=j$,  is the polynomial 
$\MSylv^{k,\ell}(P,Q)(\BU)$, where $\BU$ is a block of indeterminates of cardinality $p-j$,
\begin{align}
\MSylv^{k,\ell}(P,Q)(\BU)=
 \sum_{\genfrac{}{}{0pt}{}{\BK \subset_k \BP}
      {\BL \subset_\ell \BQ}} s_{\BK} s_{\BL} \frac{   V[(\BQ \setminus \BL) \|
     (\BP \setminus \BK)]V[\BL \| \BK \| \BU)}{ V[\BP]  V[\BQ] V(\BU)}
\end{align}

\end{definition}
In particular
\begin{align}
\MSylv^{j,0}(P,Q)(\BU)=
 \sum_{\BK \subset_j \BP} s_{\BK} \frac{ V[\BQ\|(\BP \setminus \BK)]}{V[\BQ] }
    \frac{ V[\BK  \| \BU)}{V[\BP] V(\BU)}
\end{align}

The following proposition gives the relationship between multi Sylvester double sums and Sylvester double sums.
\begin{proposition}\label{lienMsylvSylv} Denoting $\BU=U\|\BU'$ with $\BU'$ a block of $p-j-1$ indeterminates,

 $\Sylv^{k,\ell}(P,Q)(U)$ is the coefficient of $\displaystyle{\prod_{U'\in \BU'}^{}U^{\prime j}}$ in $\MSylv^{k,\ell}(P,Q)(\BU)$.

\end{proposition}

The proof of Proposition \ref{lienMsylvSylv} is based on the following Lemma.

\begin{lemma}\label{coeff0}
$ V[\BK\|U)$ is the coefficient of
$\displaystyle{\prod_{U'\in \BU'} U'^k}$ in $\displaystyle{\frac{  V[\BK\|U\|\BU')}{V(U\|\BU')}}$.
\end{lemma}

\begin{proof}
\begin{align*}
\frac{ \partial^{[\BK]}V(X_\BK\|U\|\BU')}{V(U\|\BU')} &=\frac{ \partial^{[\BK]}(V(X_\BK\|U)\Pi(\BU',X_\BK)\Pi(\BU',U)V(\BU'))}{\Pi(\BU',U)V(\BU')}\\
 &= \partial^{[\BK]}(V(X_\BK\|U)\Pi(\BU',X_\BK))
\end{align*}
Noting that 
\begin{multline*}
\partial^{[\BK]}
\left(V(X_\BK\|U)\Pi(\BU',X_\BK)\right)=
\partial^{[\BK]}V(X_\BK\|U)
\times\Pi(\BU',X_\BK)
+\sum_r V_r(X_\BK,U)\Pi_r(\BU',X_\BK)
\end{multline*}
where each $\Pi_r(\BU',X_\BK)$ is obtained by partial derivation of $\Pi(\BU',X_\BK)$ with respect to at least one variable in $X_\BK$, it is clear that the degree of some  $U'\in \BU'$ in 
$\Pi_r(\BU',X_\BK)$ is less than $k$.
The claim follows, substituting $\BK$ to $X_\BK$.
\end{proof}
\begin{proof}[Proof of Proposition \ref{lienMsylvSylv}]
The coefficient of $\displaystyle{\prod_{U'\in \BU'}^{}U^{\prime j}}$ in
$\displaystyle{\frac{  V[\BL\|\BK\|U\|\BU')}{V(U\|\BU')}}$ is $V[\BL\|\BK \| U)$ by Lemma \ref{coeff0}. The coefficient of $\displaystyle{\prod_{U'\in \BU'}^{}U^{\prime j}}$ in $\MSylv^{k,\ell}(P,Q)(\BU)$ is 
$$
\sum_{\genfrac{}{}{0pt}{}{\BK \subset_k \BP}
      {\BL \subset_\ell \BQ}} s_\BK s_\BL \frac{V[(\BQ\setminus\BL) \| (\BP\setminus\BK)] V[\BL\| \BK\| U)}{V[\BP] V[\BQ]}\
=\Sylv^{k,\ell}(P,Q)(U)$$ by Proposition \ref{lemme2bis}.
\end{proof}

\subsection{Computation of (multi) Sylvester double sums for ${j\ge q}$} \label{subsec:mult}

\begin{proposition}\label{enplus}
If $q\le j < p$
$$\MSylv^{j,0}(P,Q)(\BU)=(-1)^{j(p-j)}\prod_{U\in \BU}Q(U)$$
\end{proposition}

\begin{proof}
The polynomial $\displaystyle{\prod_{U\in \BU}Q(U)}$
 is a symmetric polynomial in $\BU$ of multidegree $q,\ldots,q$, so at most $j,\ldots,j$. Its coordinates in the basis $\mathcal{B}_{\BP,j}(\BU)$ are, for $\BK\subset_j\BP$,
$(-1)^{j(p-j)}s_\BK h(\BP\setminus\BK)$ where 
 $$h(X_{\BP\setminus\BK})=\displaystyle{\partial^{[\BP\setminus\BK]}\left(V(X_{\BP\setminus\BK})\prod_{X\in X_{\BP\setminus\BK}}Q(X)\right)}$$
  by Proposition
 \ref{interpolationdata},
 and moreover
$$h(\BP\setminus\BK)=\frac {V[\BQ\|(\BP \setminus \BK) ]}{V[\BQ]}
$$ by Lemma \ref{prodpratique}.2.

So, the polynomials  $\MSylv^{j,0}(P,Q)(\BU)$ and $\displaystyle{(-1)^{j(p-j)}\prod_{U\in \BU}Q(U)}$ have the same coordinates in the basis  $\mathcal{B}_{\BP,k}(\BU)$ and are equal. 
\end{proof}

As a corollary
\begin{proposition}\label{preouf} \
\begin{enumerate}
\item  $\Sylv^{p-1,0}(P,Q)(U)=(-1)^{p-1} Q(U)$
\item  For any  $q<j<p-1$, 
$\Sylv^{j,0}(P,Q)(U)=0$
\item  
$\Sylv^{q,0}(P,Q)(U)=(-1)^{q (p-q)}  Q(U)$
\end{enumerate}
\end{proposition}
\begin{proof}~

\begin{enumerate}
\item 
 For $j=p-1$ Proposition \ref{enplus} is exactly
$$
\Sylv^{p-1,0}(P,Q)(U)=(-1)^{p-1}Q(U).
$$.
\item  If $q<j<p-1$, denoting $\BU=U\|\BU'$ with $\BU'$ a block of $p-j-1$ indeterminates,
the coefficient of $\displaystyle{\prod_{U'\in \BU'} U'^j}$ in $\displaystyle{\prod_{U'\in \BU}Q(U')}$ is equal to $0$, so $\Sylv^{j,0}(P,Q)(U)=0$ applying Proposition \ref{enplus}. 
\item  From  Proposition \ref{enplus} and Proposition \ref{lienMsylvSylv}, denoting $\BU=U\|\BU'$ with $\BU'$ a block of $p-q-1$ indeterminates, we know that $\Sylv^{q,0}(P,Q)(U)$ is equal to the coefficient  of $\displaystyle{\prod_{U'\in \BU'} U'^q}$ in $(-1)^{q(p-q)}Q(U)\prod_{U'\in \BU'}Q(U')$. This coefficient is exactly $(-1)^{q(p-q)}Q(U)$.
\qedhere
\end{enumerate}
\end{proof}

\begin{proposition}\label{jplusgrandqueq}
If $\ell \le q\leq k+\ell=j<p$ then 
$$
\MSylv^{k,\ell}(P,Q,\BU)=(-1)_{}^{\ell(p-j)}\comb q \ell \MSylv^{j,0}(P,Q,\BU)
$$
\end{proposition}
\begin{proof}
Let $\BL\subset_\ell\BQ$ and $\BU'=(U'_1,\ldots,U'_{p-k})$; the polynomial $\displaystyle{\frac{V[(\BQ\setminus\BL)\|\BU')}{V(\BU')}}$ is a symmetric polynomial in the indeterminates $\BU'$ of degree at most $q-\ell \leq k$ in each indeterminate $U'_i, 1\leq i\leq p-k$. So, we can write this polynomial in the basis 
$\mathcal{B}_{\BP,k}(\BU')$
$$
\frac{V[(\BQ\setminus\BL)\|\BU')}{V(\BU')}=\sum_{\BK\subset_k\BP}g_\BK\frac{V[\BK\|\BU')}{V[\BP]V(\BU')}
$$
where, by Proposition \ref{interpolationdata}
$$
g_\BK=
(-1)^{k(p-k)}s_\BK V[(\BQ\setminus\BL)\|(\BP\setminus\BK)].
$$
We deduce from this
$${V[(\BQ\setminus\BL)\|\BU')}=\sum_{\BK\subset_k\BP}(-1)^{k(p-k)}s_\BK{V[(\BQ\setminus\BL)\|(\BP\setminus\BK)]}\frac{V[\BK\|\BU')}{V[\BP]}
$$
We replace  $\BU'$ by $\BU'=Y_\BL\|\BU$, where $\BU$ is a set of $p-j$ indeterminates, derivate with respect to $\partial^{[\BL]}$ and replace $Y_\BL$ by $\BL$; we obtain
\begin{align*}
\label{10}
{V[(\BQ\setminus\BL)\|\BL\|\BU)}&=\sum_{\BK\subset_k\BP} (-1)^{k(p-k})s_\BK{V[(\BQ\setminus\BL)\|(\BP\setminus\BK)]}\frac{V[\BK\|\BL\|\BU)}{V[\BP] }\\&=(-1)^{k\ell} (-1)^{k(p-k)}\sum_{\BK\subset_k\BP}s_\BK{V[(\BQ\setminus\BL)\|(\BP\setminus\BK)]}\frac{V[\BL\|\BK\|\BU)}{V[\BP] }
\end{align*}
As
 $$V[(\BQ\setminus\BL)\|\BL\|\BU)=
s_\BL
V[\BQ\|\BU),$$ we have
$$
{V[\BQ\|\BU)}=\sum_{\BK\subset_k\BP}(-1)^{k(p-j)}s_\BK s_\BL{V[(\BQ\setminus\BL)\|(\BP\setminus\BK)]}\frac{V[\BL\|\BK\|\BU)}{V[\BP]}$$
and
$$
\frac{V[\BQ\|\BU)}{V(\BU)}=\sum_{\BK\subset_k\BP}(-1)_{}^{k(p-j)}s_\BK s_\BL V[(\BQ\setminus\BL)\|(\BP\setminus\BK)]\frac{V[\BL\|\BK\|\BU)}{V[\BP] V(\BU)}.$$

The polynomial $\displaystyle{\frac{V[\BQ\|\BU)}{V(\BU)}}$ is a symmetric polynomial in the indeterminates $\BU$ of degree at most 
 $q\le j$ in each indeterminate $U_i, 1\leq i\leq p-j$. So, we can write it in the basis 
 $\mathcal{B}_{\BP,j}(\BU)$
$$\frac{V[\BQ\|\BU)}{V(\BU)}=\sum_{\BW\subset_j\BP} g_\BW\frac{V[\BW\|\BU)}{V[\BP]V(\BU)}
$$
where by Proposition \ref{interpolationdata}
$$
g_\BW=
(-1)^{j(p-j)}
s_\BW V[\BQ\|(\BP\setminus\BW)].
$$
So
\[
\sum_{\BW\subset_j\BP} s_\BW V[\BQ\|(\BP\setminus\BW)]V[\BW\|\BU)=\sum_{\BK\subset_k\BP}(-1)_{}^{\ell(p-j)}s_\BK s_\BL  V[(\BQ\setminus\BL)\|(\BP\setminus\BK)]V[\BL\|\BK\|\BU).
\]
It follows 
\begin{multline*}
\sum_{\genfrac{}{}{0pt}{}{\BW \subset_j \BP}
      {\BL \subset_\ell \BQ}} s_\BW V[\BQ\|(\BP\setminus\BW)] V[\BW\|\BU)=(-1)_{}^{\ell(p-j)}\sum_{\genfrac{}{}{0pt}{}{\BK \subset_k \BP}
      {\BL \subset_\ell \BQ}} s_\BK s_\BL V[(\BQ\setminus\BL)\|(\BP\setminus\BK)] V[\BL\|\BK\|\BU)
\end{multline*} 
\begin{multline*}
\comb q \ell \sum_{\BW\subset_j\BP} s_\BW V[\BQ\|(\BP\setminus\BW)] V[\BW\|\BU)=(-1)_{}^{\ell(p-j)}\sum_{\genfrac{}{}{0pt}{}{\BK \subset_k \BP}
      {\BL \subset_\ell \BQ}}s_\BK s_\BL V[(\BQ\setminus\BL)\|(\BP\setminus\BK)] V[\BL\|\BK\|\BU)
\end{multline*}
\begin{multline*}
\comb q \ell \frac{\sum_{\BW\subset_j\BP} s_\BW V[\BQ\|(\BP\setminus\BW)]V[\BW\|\BU)}{V[\BP]V[\BQ]V(\BU)}=(-1)_{}^{\ell(p-j)}\sum_{\genfrac{}{}{0pt}{}{\BK \subset_k \BP}
      {\BL \subset_\ell \BQ}}s_\BK s_\BL\frac{V[(\BQ\setminus\BL)\|(\BP\setminus\BK)]V[\BL\|\BK\|\BU)}{V[\BP] V[\BQ] V(\BU)}
\end{multline*}
and
\begin{align*}
 \MSylv^{k,\ell}(P,Q,\BU)=(-1)_{}^{\ell(p-j)}\comb q \ell \MSylv^{j,0}(P,Q,\BU) .&\qedhere
\end{align*}
\end{proof}

\begin{corollary}\label{lienentreSylv}For $q\leq j<p$,
\[\Sylv^{k,\ell}(P,Q)=(-1)^{\ell(p-j)}\comb q \ell \Sylv _{}^{j,0}(P,Q)\]
\end{corollary}
\begin{proof}
Immediate using Proposition \ref{jplusgrandqueq} and Proposition \ref{lienMsylvSylv}.
\end{proof}
\begin{proposition}\label{ouf} \
\begin{enumerate}
\item  
For any $(k,\ell)$ with $q=k+\ell$,
$$\Sylv^{k,\ell}(P,Q)(U)=(-1)^{k(p-q)} \comb q k Q$$
\item  For any $(k,\ell)$ with $\ell \le q,j=k+\ell$ with $q<j<p-1$, 
$$\Sylv^{k,\ell}(P,Q)(U)=0$$
\item For any  $(k,\ell)$ with $\ell \le q,k+\ell=p-1$, $$\Sylv^{k,\ell}(P,Q)(U)=(-1)^{k}\comb q \ell Q(U)$$
\end{enumerate}
\end{proposition}
\begin{proof}~
Follows from Corollary \ref{lienentreSylv} and Proposition \ref{preouf}.
\qedhere
\end{proof}

\subsection{Fundamental property of (multi) Sylvester double sums}
\label{subsec:fonda}

This section is essentially devoted to the proof of Theorem \ref{theo4}, which is a fundamental
property of Sylvester double sums: up to a constant Sylvester double sums $\Sylv^{r,j-r}(P,Q)$ depend only on $j$. Such a result has been already given for $q\leq j<p$ by Corollary \ref{lienentreSylv}.

\begin{theorem}\label{theo4} 
 If  $k \in \mathbb{N}$, $\ell \in \mathbb{N}$, $k+\ell=j<q<p$
$$
\Sylv^{k,\ell} (P,Q)(U)  =( - 1 )^{\ell(p-j)}\comb j \ell\Sylv^{j,0}(P, Q)(U).
$$
\end{theorem}

We, in fact, prove  Theorem \ref{theo4} as a corollary of a multivariate version (Theorem \ref{theo4mult}).
The proof of Theorem \ref{theo4mult} uses in an essential way the Exchange Lemma coming from \cite{KSV}.

\begin{theorem}\label{theo4mult} 
 If  $k \in \mathbb{N}$, $\ell \in \mathbb{N}$, $k+\ell=j<q<p$, and $\BU$ a set of $p-j$ indeterminates, 
$$
\MSylv^{k,\ell} (P,Q)(\BU)  =( - 1 )^{\ell(p-j)}\comb j \ell \MSylv^{j,0}(P, Q)(\BU).
$$
\end{theorem}

To prove Theorem \ref{theo4mult} for $j<q$, we need a lemma
\begin{lemma}\label{coeff}
Let $\BK \subset_k \BP, \BL \subset_\ell \BQ$ and $\BU=U_1,\ldots,U_u$ an ordered set of  variables.
Then $(-1)^{u(u-1)/2}V[\BL\|\BK]$ is the coefficient of the leading monomial  $$\prod_{i=1}^u U_i^{k+\ell+u-i}$$ of 
$V[\BL\|\BK\|\BU)$ with respect to the lexicographical ordering.
\end{lemma}
\begin{proof}
\begin{align*}
V[\BL\|\BK\|\BU)&=\partial^{[\BK]}\partial^{[\BL]}(V(Y_\BL\|X_\BK\|\BU))(\BK,\BL,\BU)\\
                               &=V(\BU)\partial^{[\BK]}\partial^{[\BL]}(V(Y_\BL\|X_\BK)\Pi(\BU,Y_\BL\|X_\BK))(\BK,\BL)
.\end{align*}The coefficient of $\displaystyle{\prod_{i=1}^u U_i^{k+\ell+u-i}}$ in $V[\BL\|\BK\|\BU)$ is $(-1)^{u(u-1)/2}$ multiplied by the coefficient of $\displaystyle{\prod_{i=1}^u U_i^{k+\ell}}$ in $\partial^{[\BK]}\partial^{[\BL]}(V(Y_\BL\|X_\BK)\Pi(\BU,X_\BK\|Y_\BL)(\BK,\BL)$. This coefficient is $$\partial^{[\BK]}\partial^{[\BL]}V(Y_\BL\|X_\BK)(\BK,\BL)=V[\BL\|\BK];$$ indeed if any derivation is done on $\Pi(\BU,X_\BK\|Y_\BL)$, with respect to $\BK$ or $\BL$, the degree in at least one indeterminate $U_i\in\BU$ decreases strictly.
\end{proof}

\begin{proof}[Proof of Theorem \ref{theo4mult} ] We have $j<q$.
Let $\BU'$ be a block of $p-\ell$ indeterminates. On one hand, 
$$
\begin{array}{rcl}
\displaystyle{\sum_{\BT\subset_\ell\BP }\Pi(X_{\BP\setminus \BT},Y_\BQ)\frac{\Pi(\BU',X_\BT)}{\Pi(X_{\BP\setminus \BT},X_\BT)}}&=&\displaystyle{\sum_{\BT\subset_\ell\BP }\frac{V(Y_\BQ\|X_{\BP\setminus \BT})V(X_\BT\|\BU')}{V(Y_\BQ)V(\BU')V(X_\BT\|X_{\BP\setminus \BT})}}\\
&=&\displaystyle{(-1)^{\ell(p-\ell)}\sum_{\BT\subset_\ell\BP }s_\BT\frac{V(Y_\BQ\|X_{\BP\setminus \BT} ) V(X_\BT\|\BU')}{V(Y_\BQ)V(\BU')V(X_\BP)}}.
\end{array}
$$

On the other hand,
$$\begin{array}{rcl}
\displaystyle{\sum_{\BL\subset_\ell\BQ} \Pi(X_\BP,Y_{\BQ\setminus \BL})\frac{\Pi(\BU',Y_\BL)}{\Pi(Y_\BL,Y_{\BQ\setminus \BL})}}
&=&
\displaystyle{\sum_{\BL\subset_\ell\BQ}\frac{V(Y_{\BQ\setminus \BL}\|X_\BP)V(Y_\BL\|\BU')}{V(X_\BP)V(\BU')V(Y_{\BQ\setminus \BL}\|Y_\BL)}}\\
&=&
\displaystyle{\sum_{\BL\subset_\ell\BQ}s_\BL\frac{V(Y_{\BQ\setminus \BL}\|X_\BP)V(Y_\BL\|\BU')}{V(X_\BP)V(\BU')V(Y_\BQ)}}
\end{array}$$
From the Exchange Lemma in \cite{KSV}, we can write 
\begin{align}\label{equation1}
\sum_{\BT\subset_\ell\BP }\Pi(X_{\BP\setminus \BT},Y_\BQ)\frac{\Pi(\BU',X_\BT)}{\Pi(X_{\BP\setminus \BT},X_\BT)}=\sum_{\BL\subset_\ell\BQ} \Pi(X_\BP,Y_{\BQ\setminus \BL})\frac{\Pi(\BU',Y_\BL)}{\Pi(Y_\BL,Y_{\BQ\setminus \BL})}
\end{align}
So, we deduce from (\ref{equation1})
\begin{align}
\sum_{\BT\subset_\ell\BP }s_\BT\frac{V(Y_\BQ\|X_{\BP\setminus \BT})V(X_\BT\|\BU')}{V(Y_\BQ)V(\BU')V(X_\BP)}=
(-1)^{\ell(p-\ell)}\sum_{\BL\subset_\ell\BQ}s_\BL\frac{V(Y_{\BQ\setminus \BL}\|X_\BP)V(Y_\BL\|\BU')}{V(X_\BP)V(\BU')V(Y_\BQ)}
\end{align}
 \begin{align} \label{depart}
\sum_{\BT\subset_\ell \BP}s_\BT V(Y_\BQ\|X_{\BP\setminus \BT})
V(X_\BT \|\BU') =(-1)^{\ell(p-\ell)}\sum_{\BL\subset_\ell \BQ}s_\BL V(  Y_{\BQ\setminus \BL}\|X_\BP) V(Y_\BL\|\BU')
\end{align}
Hence, derivating with respect to $\BQ$ and substituting $\BQ$ to $Y_\BQ$,
\begin{align} \label{depart1}
\sum_{\BT\subset_\ell \BP}s_\BT V[\BQ\|X_{\BP\setminus \BT})
V(X_\BT \|\BU') =(-1)^{\ell(p-\ell)}\sum_{\BL\subset_\ell \BQ}s_\BL V[(\BQ\setminus \BL) \|X_\BP)V[\BL\|\BU')
\end{align}

 We fix $\BK\subset_k\BP$.
 The total degree  with respect to $X_\BK$ of 
 $V[\BQ\setminus \BL\|X_\BP  ] V[\BL\|\BU')$ is $$d_1=\comb k 2+k (p+q-j).$$
  Denoting, for any $\BT\subset_\ell\BP$, $c$ the cardinality of $\BK\cap \BT$, we note that
the cardinality of $(\BP\setminus \BT)\cap \BK$ is $k-c$.

 So, the total degree  with respect to $X_\BK$ of 
 $V[\BQ\|X_{\BP\setminus \BT})
V(X_\BT \|\BU')$ is 
$$d_{2,c}=\comb {k-c} 2 +(k-c)(p+q-j+c)+\comb c 2+c(p-c),$$ i.e.
$$d_{2,c}=\comb k 2+k(p+q-j)-c(q-j+c)$$ and 
$$
d_1-d_{2,c}^{}=c(q-j+c)
$$
So $d_{2,c}<d_1$ if $c>0$ and $d_{2,c}=d_1$ if $c=0$ i.e. if $\BT  \subset \BP \setminus \BK$. This implies that subsets $\BT$ which intersect $\BK$ don't contribute to the  homogeneous part of total degree $d_1$ in $X_\BK$ on the left side of (\ref{depart1}).

 Note that, if $\BT\subset_\ell\BP\setminus\BK$, 

$$V[\BQ\|X_{\BP \setminus \BT})=r_{\BK,\BT} V[\BQ\|X_{\BP \setminus (\BK\cup \BT)}\| X_{\BK}),$$ where $r_{\BK,\BT}$ is the signature of the permutation $\rho_{\BK,\BT}$ taking the ordered set ${\BP}\setminus {\BT}$ to the ordered set $(\BP\setminus ({\BK}\cup {\BT}))\| {\BK}$.

\medskip We can also write
$$V[(\BQ\setminus\BL)\|X_\BP)=s_\BK V[(\BQ\setminus\BL)\|X_{\BP\setminus \BK}\|X_\BK)$$

\medskip If $X_\BK=X_{u_1},\ldots, X_{u_k}$, taking the coefficient of $\prod_{i=1}^k X_{u_k}^{k-i+p+q-j}$
in both sides of (\ref{depart1}) gives, by Lemma \ref{coeff}
\begin{multline} \label{suite}
\sum_{\BT\subset_\ell \BP\setminus \BK}r_{\BK,\BT}s_\BT
 V[\BQ\|X_{\BP\setminus (\BK\cup \BT)} )
V(X_\BT \|\BU') =(-1)^{\ell(p-\ell)}\sum_{\BL\subset_\ell \BQ}s_\BK s_\BL V[(\BQ\setminus \BL)\|X_{\BP\setminus \BK} ) V[\BL\|\BU').
\end{multline}
Derivating both sides of (\ref{suite}) with respect to $\partial^{[\BP \setminus \BK]}$ and replacing $X_{\BP\setminus \BK}$ by $\BP \setminus \BK$, followed by replacing $\BU'$ by $X_\BK\| \BU$, where $\BU$ is a set of $p-j$ indeterminates, derivating with respect to 
 $\partial^{[\BK]}$ and replacing $X_\BK$ by $\BK$ gives
 \begin{multline} \label{suitebis}
\sum_{\BT\subset_\ell \BP\setminus \BK}r_{\BK,\BT}s_\BT V[\BQ\|(\BP\setminus (\BK\cup \BT))]
V[\BT \|\BK\|\BU) =(-1)^{\ell(p-\ell)}\sum_{\BL\subset_\ell \BQ}s_\BK s_\BL V[(\BQ\setminus \BL)\|(\BP\setminus \BK)] V[\BL\|\BK\|\BU).
\end{multline}
Summing with respect to $\BK$, ve get
$$
\sum_{\genfrac{}{}{0pt}{}{\BK \subset_k \BP}
      {\BT\subset_\ell \BP\setminus \BK}}
   r_{\BK,\BT}s_\BT V[\BQ\|(\BP\setminus (\BK\cup\BT))]
V[\BT \|\BK\|\BU) =(-1)^{\ell(p-\ell)}\MSylv^{k,\ell}(P,Q)(U)V[\BP] V[\BQ]V(\BU).
$$
Denote $\BW$ the set $\BK\cup\BT$ ordered by the induced order on $\BP$.
Let $\tau_{\BK,\BT}$ be the permutation sending the ordered set $\BP\setminus \BW \|\BW$ to the  ordered set $\BP\setminus \BW \| \BT\|\BK$ and $t_{\BK,\BT}$ its signature.
We deduce 
$$
\sum_{\genfrac{}{}{0pt}{}{\BK \subset_k \BP}
      {\BT\subset_\ell \BP\setminus \BK}}
t_{\BK,\BT}r_{\BK,\BT}s_\BT V[\BQ\|(\BP\setminus \BW) ]
V[\BW\|\BU)=
(-1)^{\ell(p-\ell)}\MSylv^{k,\ell}(P,Q)(U)V[\BP] V[\BQ]V(\BU).
$$

We remark that $t_{\BK,\BT}r_{\BK,\BT}s_{\BT}=(-1)^{k\ell}s_\BW$.
Indeed,
denoting  $\iota_{\BK,\BT}$  the permutation sending the ordered set $(\BP\setminus \BW) \| \BT\|\BK$ to the  ordered set $(\BP\setminus \BW) \| \BK\|\BT$, with signature  $(-1)^{k\ell}$,  and  by $\rho'_{\BK,\BT}$ the permutation sending the ordered set $(\BP\setminus\BT)\|\BT$ to the ordered set $(\BP\setminus(\BK\cup \BT))\|\BK\|\BT$, with signature   $r_{\BK,\BT}$,
we have the follwing sequence of permutations
$$\begin{array}{rrcl}
\sigma_\BW:& \BP&\longleftrightarrow &(\BP\setminus\BW)\|\BW\\
\tau_{\BK,\BT}	:& (\BP\setminus\BW)\|\BW&\longleftrightarrow &(\BP\setminus\BW)\|\BT\|\BK\\
\iota_{\BK,\BT}	:& (\BP\setminus\BW)\|\BT\|\BK&\longleftrightarrow &(\BP\setminus\BW)\|\BK\|\BT\\
{\rho'}_{\BK,\BT}^{-1}:&(\BP\setminus\BW)\|\BK\|\BT&\longleftrightarrow&(\BP\setminus\BT)\|\BT\\
\sigma_\BT^{-1}:&(\BP\setminus\BT)\|\BT&\longleftrightarrow&\BP\quad\quad,
\end{array}$$
with $\sigma_\BT \circ {\rho'}_{\BK,\BT}^{-1}\circ \iota_{\BK,\BT} \circ \tau_{\BK,\BT} \circ \sigma_\BW^{-1}={\rm Id}$.

 Noting that there are $\comb j \ell$ ways of decomposing
$\BW\subset_j \BP$ as $\BW=\BK\cup \BT$, we get
 $$
\MSylv^{k,\ell} (P,Q)(\BU)  = ( - 1 )^{\ell(p-j)} \comb j \ell\MSylv^{j,0}(P, Q)(\BU).\qedhere 
$$
 \end{proof}

\begin{proof}[Proof of Theorem \ref{theo4}]
Theorem \ref{theo4} is an immediate consequence of Theorem \ref{theo4mult}, by applying Proposition \ref{lienMsylvSylv}.\end{proof}

\section{Sylvester double sums and remainders}\label{sec:remainders}

In Section \ref{sec:remainders} we give a relationship between the Sylvester double sums of $P,Q$ and those of $Q,R$ where $R$ is the opposite of the remainder of $P$ by $Q$ in the Euclidean division.

We are now dealing with not necessarily monic polynomials.

\begin{definition}\label{defnonmonic}
Let $P$ be a polynomial of degree  $p$ 
which leading coefficient is denoted  $\lc(P)$.
Let $Q$ be a polynomial of degree  $q$ which leading coefficient is denoted $\lc(Q)$.

\noindent Let $(k,\ell)$ with $j=k+\ell\le p$
be a pair of natural numbers. We define
$$\Sylv^{k,\ell}(P,Q)(U)=
{\lc(P)^{q-j}\lc(Q)^{p-j}}\Sylv^{k,\ell}\left(\frac P {\lc(P)},\frac Q{\lc(Q)}\right)(U)
$$
\end{definition}

\begin{remark}{}\label{proprietes}
{\rm Note that if $k \in \mathbb{N}$, $\ell \in \mathbb{N}$,$\ell \le q$, $k+\ell=j<q$
$$\Sylv^{k,\ell} (P,Q)(U)  = \left( - 1 \right)^{\ell(p-j)} \comb j \ell\Sylv^{j,0}(P, Q)(U)$$
follows immediately from Theorem \ref{theo4} and Definition \ref{defnonmonic}.
}
\end{remark}

We use Notation \ref{notmultiset} to define the ordered sets $\BP$ and $\BQ$ representing the multisets of roots 
of $P$ and $Q$.

Rewriting Lemma \ref{prodpratique} in the non monic case, we get  Lemma \ref{prodpratiquebis}.

\begin{lemma}\label{prodpratiquebis}\

\begin{enumerate}
\item  For $\BL\subset_\ell\BQ$,
  defining 
$$f(Y_{\BQ\setminus \BL})=(-1)^{p(q-\ell)}\partial^{{[\BQ\setminus\BL]}}\left(V(Y_{\BQ\setminus\BL})\prod_{Y \in Y_{\BQ\setminus\BL}}P(Y)\right),$$
we have
$$
f(\BQ\setminus\BL)= \lc(P)^{q-\ell}\frac{V[(\BQ\setminus\BL)\|\BP)]}{V[\BP]}
$$
\item For $\BK\subset_k\BP$,
 defining
$$g( X_{\BP\setminus\BK})=\partial^{{[\BP\setminus\BK]}}\left(V(X_{\BP\setminus\BK})\prod_{X\in X_{\BP\setminus\BK}}Q(X) \right),$$ 
we have
$$
g{(\BP\setminus\BK)}= \lc(Q)^{p-k}\frac{V[\BQ\|(\BP\setminus\BK)]}{V[\BQ]}
$$
\end{enumerate}
\end{lemma}

Similarly, reewriting  Proposition \ref{preouf} in the non monic case, we get  Proposition \ref{preoufbis}.

\begin{proposition}\label{preoufbis} \
\begin{enumerate}
\item  $\Sylv^{p-1,0}(P,Q)(U)=(-1)^{p-1}\lc(P)^{q-p+1} Q(U)$
\item  For any  $q<j<p-1$, 
$\Sylv^{j,0}(P,Q)(U)=0$
\item  
$\Sylv^{q,0}(P,Q)(U)=(-1)^{q (p-q)} \lc(Q)^{p-q-1} Q(U)$
\end{enumerate}
\end{proposition}

We proceed now to the proof of Proposition \ref{prorecurrence} wich is the main result of Section \ref{sec:remainders}.

\begin{proposition}\label{prorecurrence}
Let $R=-\rm{Rem}(P,Q).$
If  $j\in \mathbb{N}$,  $j< q$
\begin{itemize}
\item If $R=0$, $\Sylv^{j,0}(P,Q)(U) = 0$.
\item If $R\not=0$,
$\Sylv^{j,0}(P,Q)(U) = (-1)^{q(p-q)}\lc(Q)^{p-r}\Sylv^{j,0}(Q,R)(U)$
\end{itemize}
\end{proposition}

\medskip
The following elementary lemma plays a key role in the proof of Proposition
\ref{prorecurrence}.

\begin{lemma}\label{remarque}
Let $R=-\rem(P,Q).$
For every $y_{i,j}\in \BQ$, $0\leq j'< j,$ 
$$P^{[j']}(y_i)=-R^{[j']}(y_i)$$
\end{lemma}
\begin{proof}
Write $P=CQ-R$, derivate $j'$ times and evaluate at $y_i$.
\end{proof}

\begin{proof}[Proof of Proposition \ref{prorecurrence}]~
If $R=0$,  $\displaystyle{\Sylv^{0,j}\left(\frac{P}{\lc(P)},\frac{Q}{\lc(Q)}\right)(U)=0}$ follows from Corollary \ref{Regal0}. So, $$\Sylv^{j,0}(P,Q)(U)=(-1)^{j(p-j)}\Sylv^{0,j}(P,Q)(U)=0$$.

If $R\not=0$,
 let  $r$  be the degree of $R$. Let $(z_1,\ldots,z_v)$  be an ordered set of the distinct roots of $R$ in an algebraic closure $\C$ of $\mathbb{K}$, with $z_i$ of multiplicity $\xi_i$, and, as in Notation \ref{notmultiset}, let $\BR$ be the  multiset of roots of $R$, represented by the ordered set
$$\BR=(z_{1,0},\ldots,z_{1,\xi_1-1},\ldots,z_{v,0},\ldots,z_{v,\xi_v-1}), $$
with $z_{i,j}=(z_i,j)$ for $0\leq j\leq \xi_i-1$, $\sum_{i=1}^v \xi_i=r.$

If $j\le q$, define for $\BL\subset_j \BQ$
$$
\begin{array}{rcl}
f(Y_{\BQ\setminus\BL})&=&(-1)^{p(q-j)}\partial^{[\BQ\setminus \BL]}\left( V(Y_{\BQ\setminus\BL})\prod_{Y\in Y_{\BQ\setminus\BL}}P(Y) \right)\\
h(Y_{\BQ\setminus\BL})&=&\partial^{[\BQ\setminus \BL]}\left(V(Y_{\BQ\setminus\BL})\prod_{Y\in Y_{\BQ\setminus\BL}}R(Y) \right)
\end{array}
$$
Note that $$
f(\BQ\setminus\BL)=(-1)^{(p+1)(q-j)}h(\BQ\setminus\BL)
$$
from Lemma \ref{remarque}.

So
\begin{align*}
\Sylv^{0,j}(P,Q)(U)&=
\displaystyle{\frac{\lc (P)^{q-j}\lc (Q)^{p-j}}{V[\BP]V[\BQ]}\sum_{\BL \subset_j\BQ}s_\BL V[(\BQ\setminus\BL)\|\BP]V[\BL\|U)}\\
&=\displaystyle{\frac{\lc (Q)^{p-j}}{V[\BQ]}\sum_{\BL \subset_j\BQ}s_\BL f(\BQ\setminus\BL)V[\BL\|U)}\hfill\mbox{ applying Lemma \ref{prodpratiquebis}.1}\\
&=\displaystyle{(-1)^{(p+1)(q-j)}\frac{\lc (Q)^{p-j}}{V[\BQ]}\sum_{\BL \subset_j\BQ}s_\BL h(\BQ\setminus\BL)V[\BL\|U)}\\
&=\displaystyle{(-1)^{(p+1)(q-j)}\frac{\lc (Q)^{p-j}}{V[\BQ]}\sum_{\BL \subset_j\BQ}s_\BL \frac{\lc (R)^{q-j}V[\BR\|(\BQ\setminus\BL)]}{V[\BR]}V[\BL\|U)\hfill\mbox{ applying Lemma \ref{prodpratiquebis}.2}}\\
&=\displaystyle{(-1)^{(p+1)(q-j)}\frac{\lc (Q)^{p-j}\lc (R)^{q-j}}{V[\BQ]V[\BR]}\sum_{\BL \subset_j\BQ}s_\BL V[\BR\|(\BQ\setminus\BL)]V[\BL\|U)}\\
&=\displaystyle{(-1)^{(p+1)(q-j)}\frac{\lc (Q)^{p-r}\lc (Q)^{r-j}\lc (R)^{q-j}}{V[\BQ]V[\BR]}\sum_{\BL\subset_j\BQ}s_\BL V[\BR\|(\BQ\setminus\BL)]V[\BL\|U)}\\
&=\displaystyle{(-1)^{(p+1)(q-j)}{\lc (Q)^{p-r}}\Sylv^{j,0}(Q,R)}
\end{align*}
The claim follows since 
$$\Sylv^{j,0}(P,Q)(U)=(-1)^{j(p-j)}\Sylv^{0,j}(P,Q)(U)$$ by Theorem \ref{theo4}
and $$(-1)^{j(p-j)}(-1)^{(p+1)(q-j)}=(-1)^{q(p-q)}.$$

\end{proof}

\section{Sylvester double sums and subresultants}\label{sec:subresdblsum}

 Finally we prove in this section  that Sylvester double sums coincide (up to a constant) with subresultants, by an induction on the length of the remainder sequence of $P$ and $Q$. 
 
This section is devoted to the proof of the link between double sums and subresultants which is known in the simple case (see \cite{ALPP,AHKS,RS}). 

\medskip {\bf Notation.} $\varepsilon_k=(-1)^{k(k-1)/2}$.
The sign $\varepsilon_k$ is the signature of the permutation reversing the order 
i.e. sending $1,2,\ldots,k-1,k$ to
$k, k-1,\ldots,2,1$. We have also $\varepsilon_k=1$ if $k\equiv 0,1 \mod 4$, $\varepsilon_k=-1$ if $k\equiv 2,3 \mod 4$.
As a consequence
\begin{equation} \label{Eqvarepsilon}
\varepsilon_{i+1}=(-1)^{i}\varepsilon_{i}.
\end{equation}
We have also
\begin{equation} \label{Eqvarepsilonbis}
\varepsilon_{i+j}=(-1)^{ij}\varepsilon_{i}\varepsilon_{j},
\end{equation}
which follows from the fact that
reversing $i+j$ numbers can be done in three steps: reversing the  first $i$ ones, then the  last $j$ one and placing the last $j$ numbers in front of the $i$ first.

\noindent The main theorem of this section is the following.

\begin{theorem}\label{theoreme2}
Let  $k \in \mathbb{N}$, $\ell \in \mathbb{N}$, $\ell \le q$, $k+\ell=j<p-1$
$$\Sylv^{k,\ell}(P,Q)(U)=(-1)^{k(p-j)}\varepsilon_{p-j}\comb j k\Sres_{j}(P,Q)(U).$$
\end{theorem}

\begin{remark}{}
{\rm When $j=p-1$, $\Sres_{p-1}^{}(P,Q)(U)=Q(U)$ by convention; so, as $$\displaystyle{\Sylv^{k,\ell}(P,Q)(U)=(-1)^k\comb q \ell \lc(Q)^{q-p+1} Q(U)}\mbox{ for }k+\ell=p-1,\mbox{ we get }\hfill$$
$$
\Sylv^{k,\ell}(P,Q)(U)=(-1)^k\comb q \ell \lc(Q)^{q-p+1}  \Sres_{p-1}^{}(P,Q)(U).
$$}
\end{remark}

In order to prove Theorem \ref{theoreme2}, we use an induction on the length of the remainder sequence of $P$ and $Q$ based on Proposition \ref{prorecurrence}.

Before proving Theorem \ref{theoreme2}, we recall the following properties of subresultants.
\begin{lemma}\label{rappel}~
Let $R=-\rem (P,Q).$
 $$
\begin{array}{lll}
1.\,q<j<p-1&\quad\Sres_j(P,Q)(U)&=0\\
2.\, j=q&\quad\Sres_q(P,Q)(U)&=\varepsilon_{p-q}\lc(Q)^{p-q-1}Q(U)\\
3.\, j=q-1&\quad\Sres_{q-1}(P,Q)(U)&=\varepsilon_{p-q}\lc(Q)^{p-q+1}R(U)\\
4.\, j<q-1, R\not=0&\quad\Sres_j(P,Q)(U)&=\varepsilon_{p-q}\lc(Q)^{p-r}\Sres_j(Q,R)(U)\\
5.\, j<q-1, R=0 &\quad\Sres_j(P,Q)(U)&=0
\end{array}
$$
\end{lemma}
\begin{proof}
All items follow from \cite{BPR} except the computation of $\Sres_{q-1}(P,Q)(U)$.
$\Sres_{q-1}(P,Q)(U)$ is clearly equal to $\varepsilon_{p-q+2}\lc(Q)^{p-q+1}(-R(U))$ by replacing the row of $P$ by a row of $-R$ in the Sylvester-Habicht matrix, and reversing the order of its $p-q+2$ rows.
Notice now that $\varepsilon_{p-q+2}=-\varepsilon_{p-q}$.
\end{proof}

We also recall the following similar properties of Sylvester double sums.
\begin{lemma}\label{rappelbis}
Let $R=-\rem(P,Q).$
Let  
$j\in \mathbb{N}$, $j<p-1$
$$\begin{array}{lll}
1.\,q<j<p-1&\Sylv^{j,0}(P,Q)(U)&=0,\\
2.\,j=q&\Sylv^{q,0}(P,Q)(U)&=(-1)^{q(p-q)} \lc(Q)^{p-q-1}Q(U)\\
3.\,j=q-1&\Sylv^{q-1,0} \left( P, Q \right)(U)  &=(-1)^{(q-1)(p-q+1)+p-q} \lc(Q)^{p-q+1}R(U)\\
4.\,j<q-1,
       R\not=0&\Sylv^{j,0}(P,Q)(U)&=(-1)^{q(p-q)} \lc(Q)^{p - r}
       \Sylv^{j, 0} (Q,R)(U) 
       \\
5.\, j<q-1, 
 R=0&\Sylv^{j,0}(P,Q)(U)&=0, 
     
\end{array}$$
\end{lemma}
\begin{proof}
All items follow from  Proposition \ref{preoufbis} except for the computation of $\Sylv^{q-1,0}(P,Q)(U)$.
Using   Proposition \ref{prorecurrence} and Proposition \ref{preoufbis} for $Q,R$, we obtain
$$
\Sylv^{q-1,0}(P,Q)(U)= (-1)^{q(p-q)}\lc (Q)^{p-r}\Sylv^{q-1,0}(Q,R)(U)=(-1)^{q-1+q(p-q)}\lc (Q)^{p-r+r-q+1}R(U)
$$
It remains to remark that
$(q-1)(p-q+1)+p-q=q(p-q)+q-1$.
\end{proof}

\begin{proof}[Proof of Theorem \ref{theoreme2}]
The statement for $q\le j <p-1$ follows from Lemma \ref{rappel} 1,2, 
Lemma \ref{rappelbis} 1,2  and Theorem \ref{theo4}.

The statement for $j=q-1$ follows from Lemmma \ref{rappel} 3,
Lemma \ref{rappelbis} 3, Theorem \ref{theo4} and (\ref{Eqvarepsilon}) since $\varepsilon_{p-q+1}=(-1)^{p-q}\varepsilon_{p-q}.$

For $j<q-1$ we first prove the special case
\begin{align}\label{star}
\Sylv^{j,0}(P,Q)=(-1)^{j(p-j)}\varepsilon_{p-j}\Sres_j(P,Q)
\end{align}
The proof is by induction on the length of the remainder sequence of $P,Q$.

The basic case is when $Q$ divides $P$, i.e. $R=0$, and the claim is true by 
Lemma \ref{rappel} 5,
Lemma \ref{rappelbis} 5.

Otherwise suppose, by induction hypothesis that 
\begin{align}
\Sylv^{j,0}(Q,R)=(-1)^{j(q-j)}\varepsilon_{q-j}\Sres_j(Q,R)
\end{align}
Using Lemma \ref{rappel} 5 and 
Lemma \ref{rappelbis} 5
it remains to note that 
$$(-1)^{j(p-j)}(-1)^{j(q-j)}(-1)^{q(p-q)}=(-1)^{(q-j)(p-q)}$$ and conclude by 
(\ref{Eqvarepsilonbis}) since
$
\varepsilon_{p-j}=(-1)^{(q-j)(p-q)}\varepsilon_{p-q}\varepsilon_{q-j}.
$ 

The general case for $k,\ell$ now follows from Theorem \ref{theo4}.

\end{proof}

The authors thank the referees for their relevant remarks. Special thanks to Daniel Perrucci for a very careful rereading.

\bibliography{references}

\end{document}